\documentclass[11pt]{amsart}
\usepackage{enumerate,amssymb,amsfonts,latexsym,psfrag,psfig}
\usepackage{amscd,amsmath}
\usepackage[dvips]{graphicx}
\usepackage{fullpage} 

\newcommand{\barrow}{\psfig{figure=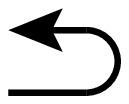,height=.10in}}
\newcommand{\weaks}{weak${\,\ast}$\ }

\newcommand{\card}{\operatorname{card}}
\newcommand{\supp}{\operatorname{supp}}
\newcommand{\diam}{\operatorname{diam}}

\def\veps{\varepsilon}
\def\Z{\mathbb{Z}}
\def\R{\mathbb{R}}

\def\N{\mathbb{N}}

\def\cG{\mathcal G}
\def\cF{\mathcal F}

\def\fM{\mathfrak M}

\def\cK{\mathcal K}

\def\hf{\hat f}
\def\tf{\tilde f}
\def\hmu{\hat\mu}
\def\tnu{\tilde\nu}

\def\tlambda{\tilde\lambda}
\def\hlambda{\hat\lambda}

\def\Ld{\Lambda_\delta}
\def\d{\delta}

\def\sph{\mathbb{S}^2}
\def\torus{\mathbb{T}^2}

\newcommand{\I}{^{-1}}
\newcommand{\ol}[1]{\overline{#1}}

 \newcommand{\pichere}[2]
 {\begin{center}\includegraphics[width=#1\textwidth]{#2}\end{center}}
 \newcommand{\lab}[3]{\psfrag{#1}[#3]{$\scriptstyle{#2}$}}

\newtheorem{Theorem}{Theorem}
\newtheorem{Lemma}[Theorem]{Lemma}

\theoremstyle{definition}

\newtheorem{Definition}{Definition}

\theoremstyle{remark}

\title{Monotone quotients of surface diffeomorphisms}
\author{Andr\'e de Carvalho}
\address{Institute for Mathematical Sciences\\
SUNY\\
Stony Brook, NY 11794-3660}
\email {andre@math.sunysb.edu}
\author{Miguel Paternain}
\address{Centro de Matem\'atica\\
Facultad de Ciencias\\
Eduardo Acevedo 1139\\
Montevideo CP 11200\\
Uruguay}
\email{miguel@cmat.edu.uy}

\begin{document}

\bibliographystyle{alpha}

%\openup 1.8\jot

\begin{abstract}
A homeomorphism of a compact metric space is {\em tight} provided
every non-degenerate compact connected (not necessarily invariant)
subset carries positive entropy.  It is shown that every
$C^{1+\alpha}$ diffeomorphism of a closed surface factors to a tight
homeomorphism of a generalized cactoid (roughly, a surface with nodes)
by a semi-conjugacy whose fibers carry zero entropy.
\end{abstract}

\maketitle
\thispagestyle{empty} \def\IMSmarkvadjust{0 pt}
\def\IMSmarkhadjust{0 pt}
\def\IMSmarkhpadding{0 pt}
\def\IMSpubltext{Published in modified form:}
\def\SBIMSMark#1#2#3{
 \font\SBF=cmss10 at 10 true pt
 \font\SBI=cmssi10 at 10 true pt
 \setbox0=\hbox{\SBF \hbox to \IMSmarkhpadding{\relax}
                Stony Brook IMS Preprint \##1}
 \setbox2=\hbox to \wd0{\hfil \SBI #2}
 \setbox4=\hbox to \wd0{\hfil \SBI #3}
 \setbox6=\hbox to \wd0{\hss
             \vbox{\hsize=\wd0 \parskip=0pt \baselineskip=10 true pt
                   \copy0 \break%
                   \copy2 \break% 
                   \copy4 \break}}
 \dimen0=\ht6   \advance\dimen0 by \vsize \advance\dimen0 by 8 true pt
                \advance\dimen0 by -\pagetotal
	        \advance\dimen0 by \IMSmarkvadjust
 \dimen2=\hsize \advance\dimen2 by .25 true in
	        \advance\dimen2 by \IMSmarkhadjust

%
%   Check for publication info
%
%  \newread\jref
  \openin2=publishd.tex
  \ifeof2\setbox0=\hbox to 0pt{}
  \else 
     \setbox0=\hbox to 3.1 true in{
                \vbox to \ht6{\hsize=3 true in \parskip=0pt  \noindent  
                {\SBI \IMSpubltext}\hfil\break
                \input publishd.tex 
                \vfill}}
  \fi
  \closein2
  \ht0=0pt \dp0=0pt
 \ht6=0pt \dp6=0pt
 \setbox8=\vbox to \dimen0{\vfill \hbox to \dimen2{\copy0 \hss \copy6}}
 \ht8=0pt \dp8=0pt \wd8=0pt
 \copy8
 \message{*** Stony Brook IMS Preprint #1, #2. #3 ***}
}

\def\IMSmarkvadjust{-30pt}
\SBIMSMark{2002/04}{October 2002}{}

\section{Introduction}
In this article, techniques from point set topology developed in the
'20's and '30's and more recent dynamical systems techniques are
brought together to construct {\em tight} models for surface
diffeomorphisms. Roughly speaking, a self-map of a metric space is
tight if any non-degenerate continuum behaves chaotically under
iteration of the map (i.e. carries positive topological entropy).  The
tight models of surface diffeomorphisms constructed here are
homeomorphisms of singular surfaces (with possibly infinitely many
nodes).  From the point of view of topological entropy, they retain
the dynamically meaningful part of the original maps, and are
`minimal' with respect to this property.

Singular surfaces (a.k.a.\ surfaces with nodes, generalized cactoids,
etc.)\ and maps between them have been present in mathematics for a
long time. They occur naturally in Thurston's theory of surface
homeomorphisms. Pseudo-Anosov maps and the generalized pseudo-Anosov
maps introduced in~\cite{dC} are examples of tight maps. Thurston's
classification theorem for surface homeomorphisms constructs minimal
complexity models in each isotopy class of homeomorphisms of the
surface. These can be viewed as piecewise affine maps on noded Riemann
surfaces: by collapsing to points the reducing curves and the finite
order components, what remains is a noded Riemann surface and the
return map to each piece is a pseudo-Anosov map. These models usually
have strictly less dynamics than other maps in the same isotopy
class. The tight models constructed here, on the other hand, preserve
the dynamics of the original map and collapse to points maximal
connected sets which are dynamically irrelevant (in the sense of
topological entropy). Under appropriate hypotheses, the Thurston and
tight models coincide: this is the case, for example, for Axiom A
surface diffeomorphisms satisfying some additional hypotheses.

An important example of the appearance of singular surfaces occurs in
Gromov's compactness theorem for $J$-holomorphic curves in symplectic
manifolds. There is some evidence that the class of tight maps on generalized
cactoids will provide the framework for the completion of the set of
pseudo-Anosov maps and for the construction of piecewise affine models
for surface diffeomorphisms in a quite general context (including
Thurston's models as a countable sub-family). Some of the techniques
needed to prove such results are similar to those used in the
proof of Gromov's compactness.

The proofs of the main theorems of the present paper use some
technical results about the dynamics on the space of continua. This
space may be thought of as a weak version of the tangent
bundle and the topological entropy of a continuum as a weak version of the Lyapunov exponent. Although the definition of entropy of non-invariant sets was
introduced by Bowen together with his definition of topological
entropy, it does not seem to have been used systematically. The
results presented in Section~\ref{sec:app} use this definition in the 
context of the space of continua and
hold in the general setting of homeomorphisms of compact metric
spaces. They are similar in flavor to some results from smooth ergodic
theory and may have some independent interest.

\smallskip
In Section~\ref{sec:back} the preliminary definitions and the
statements of the main theorems (Theorems~\ref{thm:main}
and~\ref{thm:musc}) are given. Section~\ref{sec:example} presents an
example: the quotient of Smale's horseshoe under the zero-entropy
equivalence. It is is useful to keep it in mind while reading the
proofs of the main theorems, which are given in
Section~\ref{sec:proof}. The
proof of a technical result needed in Section~\ref{sec:proof} is
postponed to Section~\ref{sec:app}.  In
Section~\ref{sec:expansive}, Theorem~\ref{thm:expansive}, giving
sufficient conditions for maps to be tight, is proved and
Section~\ref{sec:further} contains some problems for further research.

\smallskip\noindent {\bf Acknowledgements:} Both authors wish
to thank the Institute for Mathematical Sciences of SUNYSB for the
hospitality and support while this work was in progress.

\section{Background and statement of results}
\label{sec:back}

Begin by recalling Bowen's definition of topological entropy (see
\cite{Bo,G,KH}). Let $(X,d)$ be a metric space and $f\colon X \barrow$
be a uniformly continuous homeomorphism.  Define
\[d_n(x,y)= \max \left\{d(f^i(x),f^i(y))\,\,|\,\,0\leq i < n \right\}.\]
A set $E$ is $(n,\veps)$-{\em separated} if for any two distinct points
$x,y\in E$, $d_n(x,y)>\veps$. A set $F$ $(n,\veps)$-{\em spans}
another set $K$ provided that for every $x\in K$ there exists $y\in F$
such that $d_n(x,y)\leq\veps$. Let $K \subset X$ be a compact subset
and define the following quantities: $r(n,\veps,K)$ is the minimal
cardinality of a set which $(n,\veps)$-spans $K$, $s(n,\veps,K)$ is
the maximal cardinality of an $(n,\veps)$-separated subset of $K$ and
$D(n,\varepsilon,K)$ is the minimum number of sets whose
$d_n$-diameter is smaller than $\varepsilon$ and whose union covers
$K$. With these definitions, the limits below all exist and are equal
and $h(f,K)$ is defined to be equal to all of them:
\begin{eqnarray*}
h(f,K) &=& \lim_{\varepsilon \rightarrow 0}  \limsup_{n \rightarrow
\infty} \frac{1}{n}\ln r(n,\varepsilon,K)\\
       &=& \lim_{\varepsilon \rightarrow 0}  \limsup_{n \rightarrow
\infty} \frac{1}{n}\ln s(n,\varepsilon,K)\\
       &=& \lim_{\varepsilon \rightarrow 0}  \lim_{n \rightarrow
\infty} \frac{1}{n}\ln D(n,\varepsilon,K)  
\end{eqnarray*}
It is not hard to show that $h(f,f(K))=h(f,K)$ and, provided $K$
is $f$-invariant, that $h(f,K)= h(f^{-1},K)$.
Define the {\em (topological) entropy of $f$ in
$K$} or the {\em entropy carried by $K$ under $f$} to be
$h_{\pm}(K)=\max\{h(f,K),h(f^{-1},K)\}$. The {\em topological entropy}
of $f$ is defined as $h(f)=\sup h(f,K)$, where the supremum is taken
over all compact subsets $K\subset X$.

A set is {\em non-degenerate} if it contains more than one point. A
{\em continuum} is a compact connected (subset of a) metric space.
\begin{Definition}
A homeomorphism $f\colon X\barrow$ is {\em tight} if $h_\pm(C)>0$ for
every non-degenerate continuum $C\subset X$.
\end{Definition}

\begin{figure}[htbp]
\begin{center}
\pichere{0.6}{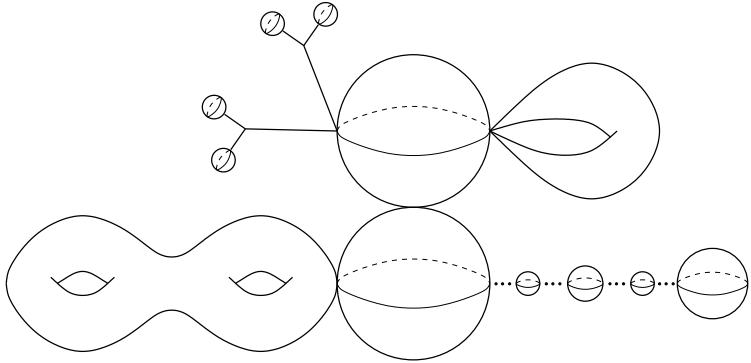}
\caption{A generalized cactoid.}
\label{fig:cactoid}
\end{center}
\end{figure}

Now {\em cactoids} and {\em generalized cactoids} are defined. This
requires several more simple definitions which are given in the list
below. Let $X$ be a connected topological space:
\begin{enumerate}[a)]
\item A point $p\in X$ for which $X\setminus \{p\}$ is not connected,
is a {\em cut point} of $X$.
\item An {\em endpoint} of $X$ is a point which has arbitrarily small
neighborhoods whose boundary is a single point.
\item A cut point $q$ {\em separates} two points $p,p'$ if it is
possible to write $X\setminus \{q\}=A\cup B$ where $p\in A, p'\in B$ and
$\ol{A}\cap B=\emptyset=A\cap\ol{B}$.
\item If $p\in X$ is neither a cut point nor an endpoint of $X$, the
set of all points which cannot be separated from $p$ by any other
point in $X$ is called a {\em (simple) link} of $X$.
\end{enumerate}
\begin{Definition}
Let $X$ be a locally connected continuum. If each simple link of $X$
is homeomorphic to the 2-sphere $\sph$ then $X$ is called a {\em
cactoid}. Let $Y$ be a space each link of which is homeomorphic to a
surface and all but finite many links are homeomorphic to $\sph$, and
$X$ be obtained by identifying finitely many pairs of points of
$Y$. Such an $X$ is called a {\em generalized cactoid}
(Figure~\ref{fig:cactoid}).
\end{Definition}

A generalized cactoid thus may have parts which are graphs or
dendrites, but all its `fat' parts have to be surfaces and, in fact,
all but finitely many have to be 2-spheres.

A map $f\colon X \barrow$ is {\em semi-conjugate} to
$g\colon Y\barrow$ if there is a continuous surjective map $\pi:X
\rightarrow Y$ such that $g \circ \pi=\pi \circ f$.  The sets of the
form $\pi^{-1}(x)$ are the {\em fibers} of the semi-conjugacy.  If
$\pi$ is a homeomorphism, then $f,g$ are {\em conjugate}.

\smallskip
The main result of this paper is  

\begin{Theorem}
\label{thm:main}
Every $C^{1+\alpha}$ diffeomorphism of a closed surface is
semi-conjugate to a tight homeomorphism of a generalized cactoid by a
semi-conjugacy whose fibers carry zero entropy.  
\end{Theorem}

This theorem is, in fact, a corollary of Theorem~\ref{thm:musc}
below. In order to state it, it is necessary to introduce the
following definition, which is central to the paper: 
\begin{Definition}
Let $f\colon X\barrow$ be a homeomorphism. Define two points $x$ and
$y$ to be {\em zero-entropy equivalent}, denoted by $x\sim y$, if there
is a continuum $C$ containing both points with $h_\pm(C)=0$. 
\end{Definition}
That this indeed defines an equivalence relation follows
from the facts that $h(f,K\cup K')=\max\{h(f,K),h(f,K')\}$ and that the
union of two connected sets with one point in common is also
connected.

It is in dimensions 1 and 2 that the zero-entropy equivalence seems to be
most interesting, usually inducing a non-trivial partition of the space.
\begin{Definition}\label{def:cactoid}
Let $X$ be a metric space and $\cG$ a partition (or {\em
decomposition}) of $X$. $\cG$ is a {\em monotone decomposition} if it
is a partition into connected sets. It is {\em upper semi-continuous}
if the sets in $\cG$ are compact and for each set $\gamma\in \cG$ and
each open set $U\supset\gamma$, there exists another open set
$V\supset\gamma$ such that every $\gamma'\in\cG$ intersecting $V$ is
contained in $U$.
\end{Definition}

The collection $\cG$ should also be thought of as the collection of
equivalence classes of the equivalence relation $x\sim y$ if and only
if $x,y\in\gamma$ for some $\gamma\in\cG$. If $X$ is a compact metric
space --- which is the case considered here --- then $\cG$ is upper
semi-continuous if and only if $x_n \sim y_n$, $x_n
\rightarrow x$ and $y_n \rightarrow y$ imply $x \sim y$.  In what
follows, we will refer interchangeably to the equivalence relation and
the corresponding decomposition.

Theorem~\ref{thm:main} is a consequence of
\begin{Theorem}
\label{thm:musc}
Let $f\colon M\barrow $ be a $C^{1+\alpha}$ diffeomorphism of a closed
surface $M$. Then the zero-entropy equivalence relation induces a
monotone upper semi-continuous decomposition of~$M$. Moreover, the
elements of this decomposition carry zero entropy. 
\end{Theorem}

The topological tool used to derive Theorem~\ref{thm:main} from
Theorem~\ref{thm:musc} is a theorem by Roberts and Steenrod~\cite{RS},
which generalizes for surfaces a result Moore~\cite{Moo} proved for
the sphere. Moore's theorem states that, if $\sim$ is a monotone upper
semi-continuous equivalence relation on the sphere $\sph$, then the
quotient space $\sph/\!\!\sim$ is a cactoid; moreover, if no
equivalence class of $\sim$ separates $\sph$, then $\sph/\!\!\sim$ is
homeomorphic to $\sph$. Roberts and Steenrod's theorem states that if
$M$ is a closed surface and $\sim$ is a monotone upper semi-continuous
equivalence relation then the quotient space is a generalized cactoid.
(The definition of generalized cactoid given here is slightly more
general than the one in~\cite{RS}, making their result simpler to state. 
In their paper, the definition does not include the possibility of
identifying finitely many pairs of points as is done in
Definition~\ref{def:cactoid}.)

\bigskip
Now a theorem giving a sufficient condition for a map to be
tight is stated.
\begin{Definition}
A homeomorphism $f\colon X\barrow$ of a metric space $X$ is {\em
continuum-expansive} if there is $\varepsilon>0$ such that, if a 
continuum $C$ satisfies $\diam f^n(C)\leq \varepsilon$ for every $n \in
\Z$, then $C$ is a point.
\end{Definition}
Expansive homeomorphisms (that is,
homeomorphisms for which there exists a constant $c>0$ such that
iteration under $f$ or $f^{-1}$ brings any two distinct points at
least $c$ apart) are clearly continuum-expansive but the converse is
not true: pseudo-Anosov homeomorphisms with 1-pronged singularities
are continuum-expansive but are not expansive (for example, if
$A\colon\torus\barrow$ a linear torus Anosov, the quotient space
$\torus/\!\!\sim$ under the identification $x\sim -x$ is homeomorphic
to the 2-sphere $\sph$ and the projection map
$\pi\colon\torus\rightarrow\sph$ is a branched covering with $4$
branch points; the map $A$ induces a map $f_A\colon \sph\barrow$ which
is a pseudo-Anosov map with 1-pronged singularities at the branch
points.)

A partial characterization of tight maps is given by
\begin{Theorem}
\label{thm:expansive}
Continuum-expansive homeomorphisms are tight. 
\end{Theorem}

\section{An example}
\label{sec:example}

Now a brief description of the zero-entropy equivalence relation for
Smale's horseshoe map (\cite{Sm}) is given. For a more detailed
discussion of this and other related examples, see~\cite{dC}.

\begin{figure}[htbp]
\begin{center}
\lab{R}{R}{l}
\lab{f}{f}{}
\lab{x}{x}{}
\lab{y}{x_0}{}
\lab{z}{x_1}{}
\pichere{0.4}{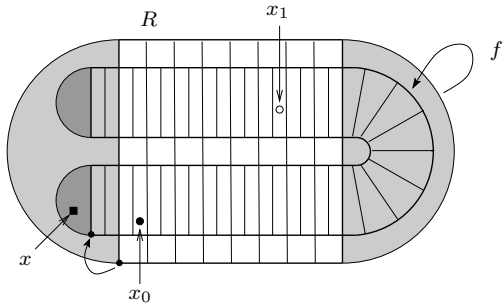}
\caption{The horseshoe map.}
\label{fig:hs}
\end{center}
\end{figure}

The horseshoe is a homeomorphism $f\colon\sph\barrow$ 
which stretches horizontally and squeezes vertically a stadium shaped
region $R$, placing it inside itself as shown in
Figure~\ref{fig:hs}. There is a repeller at
infinity whose basin contains all points outside $R$.
The horseshoe has two saddle fixed points which are
labeled $x_0$ and $x_1$ (shown as $\bullet$ and $\circ$,
respectively, in Figure~\ref{fig:hs}) and an attracting fixed point in
the shaded semi-circular region on the left denoted by $x$ (shown as
\rule{1.5mm}{1.5mm}). 

Denote by ${\mathcal H}^{u}$ and ${\mathcal H}^s$ the closures of the unstable
and stable manifolds of the fixed point~$x_0$ (or indeed of any other
periodic point, since their closures coincide) and let ${\mathcal H} =
{\mathcal H}^s\cup{\mathcal H}^u$. Equivalence classes of the zero-entropy
equivalence for the horseshoe are of four kinds:
\begin{enumerate}[a)]
\item Closures of connected components of $\sph\setminus{\mathcal H}$.
\item Closures of connected components of ${\mathcal H}^u\setminus{\mathcal H}^s$
(not already contained in sets in a)). 
\item Closures of connected components of ${\mathcal H}^s\setminus{\mathcal H}^u$
(not already contained in sets in a)).
\item Single points which are in none of the sets in a), b) or c).
\end{enumerate}

To see that these sets do not carry entropy, notice that all points in
any connected component of $\sph\setminus{\mathcal H}$ (before taking the
closure) converge to the attracting fixed point $x$. It is not hard to
see that, after taking the closure nothing more significant happens and
this shows the sets in a) indeed carry no entropy. The same holds for
sets of types b) and c). To see that any larger continuum must contain
entropy, notice that if $C$ is a connected set that contains two
distinct sets among the ones described above, then it must intersect a
Cantor set's worth of invariant manifolds, either stable or unstable
(or both). It follows that one of its $\omega$- or $\alpha$-limit sets
contains all the non-wandering set of the horseshoe and therefore one of 
$h(f,C)$ or $h(f^{-1},C)$ equals $\ln 2$. 

\begin{figure}
\pichere{0.4}{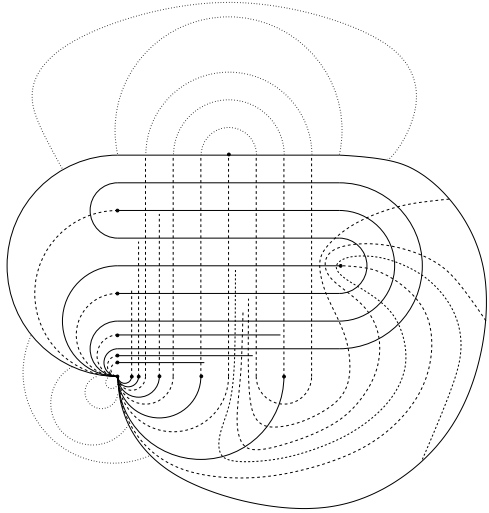}
\caption{The quotient of the sphere under the zero-entropy equivalence
relation for the horseshoe.}
\label{fig:tightHS}
\end{figure}

The quotient space is represented in Figure~\ref{fig:tightHS}. It is
again a sphere, obtained by identifying the solid boundary in the
figure along the dotted arcs from the mid-point at the top to the
corner point on the lower left. The stable and unstable manifolds of
the horseshoe project to two transverse foliations with singularities,
represented by solid and dashed lines, respectively. In fact, these
foliations carry transverse invariant measures whose product gives a
Euclidean structure on the sphere. The quotient map preserves both
foliations, dividing one of the transverse measures by 2 and
multiplying the other also by 2, so that the product measure is
invariant. This map is a {\em generalized pseudo-Anosov map}
(see~\cite{dC}).

\section{Proofs of Theorems~\ref{thm:main} and~\ref{thm:musc}}
\label{sec:proof}

The idea of the proof is roughly as follows. It is necessary to show that
carrying zero entropy is a closed condition or, alternatively, that
carrying positive entropy is an open condition: if a continuum $C$
carries (positive) entropy, then all nearby continua must also carry
entropy. Suppose then that $h_\pm(C)>0$. The first step is to obtain,
using the dynamics on the space of continua, an ergodic invariant measure with positive entropy such that arbitrarily near almost every point there is an iterate 
of $C$ with definite diameter.  
In the context of $C^{1+\alpha}$ surface diffeomorphisms, results due to Katok~\cite{Ka,KH} show 
that there must be a horseshoe near the support of 
any ergodic measure with positive entropy. This
implies that iterates of $C$ 
must eventually intersect the (un)stable set of this
horseshoe. The 2-dimensionality of the ambient space forces such
intersections to be `transversal,' so that (appropriate iterates of)
continua near $C$ must also intersect the (un)stable set of the
horseshoe and thus carry entropy.

\medskip
Begin by setting some definitions and notations.  Let $(X,d)$ be a
metric space and $\mu$ a Borel probability measure ({\em Bpm}) on
$X$. The {\em support} of $\mu$, $\supp\mu$ , is the set of points
$x\in X$ such that every neighborhood $V\ni x$ satisfies $\mu(V)>0$. A
point $x\in X$ is an {\em atom}\footnote{This definition applies to
Borel measures in separable metric spaces.} of $\mu$ if
$\mu(\{x\})\neq 0$.

The {\em Hausdorff distance} between two compact subsets $A,B\subset
X$ is
\[d_H(A,B)=\max\{\max_{a\in A}d(a,B),\max_{b\in B}d(b,A)\}\]
where, if $x\in X, C\subset X$, $d(x,C)=\inf_{c\in C}d(x,c)$. If $X$
is compact, the Hausdorff metric makes the set $\cK=\cK(X)$ of all
compact subsets of $X$ into a compact metric space.

The following lemma is one of the most important technical ingredients
in the proof of the main theorems. Its proof is given in the next
section.  Let $f\colon X\barrow$ be a homeomorphism and denote by
$\omega(C)\subset\cK$ the omega-limit set of $C\in\cK$ with respect to
the map $K\mapsto f(K)$. If $\mu$ is an $f$-invariant Bpm, $h_\mu(f)$
denotes the measure theoretic entropy of $f$ with respect to $\mu$.

\begin{Lemma}
\label{lem:basic}
Let $(X,d)$ be a compact metric space, $f\colon X\barrow$ be a
homeomorphism and $C\subset X$ be a compact set satisfying
$h(f,C)>0$. Then there exists an $f$-invariant ergodic Bpm $\mu$ on
$X$, such that $h_{\mu}(f)>0$ and such that $\mu$-almost every point
belongs to a non-degenerate set in $\omega(C)$.
\end{Lemma}

\medskip
Now the assumption of differentiability is added to the discussion
and attention is restricted to surfaces. Although the assumption that
$f$ be $C^1$ is all that is needed for some of the definitions and
results listed below, in order to use Pesin Theory --- which is
essential to our results --- it is necessary to assume some extra
regularity: {\em from now on, $M$ will denote a compact smooth
Riemannian surface without boundary and $f\colon M\barrow$ a
$C^{1+\alpha}$ diffeomorphism, with $0<\alpha\le 1$.}

\medskip
The results now described follow from Pesin
Theory~\cite{Pe,PS,Ka,KH}. The concepts needed are quite technical and
the definitions rather involved. Instead of presenting them in detail,
which would duplicate what is contained in the papers cited, only a
description of the concepts and results that will be used is given.

An $f$-invariant probability measure $\mu$ is {\em hyperbolic} if all
the Lyapunov exponents of $f$ are non-zero at $\mu$-almost every
point. Notice that if $\mu$ is an ergodic $f$-invariant Bpm on (the
surface) $M$ with $h_\mu(f)>0$ then ergodicity implies that $\mu$ has
no atoms and Ruelle's inequality implies that it is hyperbolic. This
observation will be used below.

Assume that $\mu$ is a non-atomic hyperbolic ergodic $f$-invariant
Bpm. Given $0<\delta<1$ there exists a compact {\em Pesin set}
$\Lambda_\delta$ with $\mu(\Ld)>1-\d$ and with the properties that are
now described. To avoid cluttering the notation, the dependence on
$\d$ and the point $p$ is not incorporated into it. For every
$p\in\Ld$ there exist an open neighborhood $N\ni p$, a compact
sub-neighborhood $V\ni p$, and a diffeomorphism $\Psi\colon
(-1,1)^2\to N$, with $\Psi(0,0)=p$ and $\Psi([-1/10,1/10]^2)=V$, such
that the $N$-local unstable manifolds $W^u_N(y)$ of all points $y$ in
$\Lambda_\d\cap V$ are the images under $\Psi$ of graphs of the form
$\{(v,\varphi(v))\,|\, v\in (-1,1)\}$ with small Lipschitz
constant. Any two such local unstable manifolds are either disjoint or
equal and they depend continuously on the point $y\in\Ld\cap
V$. Similarly, the $N$-local stable manifolds $W^s_N(y)$ of points $y$ in
$\Lambda_\d\cap V$ are the images under $\Psi$ of graphs of the form
$\{(\varphi(v),v)\,|\, v\in (-1,1)\}$ with small Lipschitz
constant. Any two such local unstable manifolds are either disjoint or
equal and they depend continuously on the point $y\in\Ld\cap V$.
Stable and unstable manifolds of points in $V$ having the properties
just described are called $N$-{\em admissible}.

It follows that there is a continuous product structure: given any
$x,y\in \Lambda_\d\cap V$, the intersection $W^u_N(x)\cap W^s_N(y)$ is
transversal and consists of exactly one point, which is denoted by
$[x,y]$. Define maps $\pi^s_x\colon \Lambda_\d\cap V\to W^s_N(x)$ and
$\pi^u_x\colon \Lambda_\d\cap V\to W^u_N(x)$ by $\pi^s_x(y)=[y,x]$ and
$\pi^u_x(y)=[x,y]$. Observe, however, that $[x,y]$ need not be in
$\Ld$.

Katok's Closing Lemma (the Main Lemma in~\cite{Ka}) will also be
needed: it states that arbitrarily near any recurrent point in
$\Ld\cap V$ there exist hyperbolic periodic saddles whose invariant
manifolds are $N$-admissible.

Let $A$ denote the subset of $\Ld\cap V$ consisting of all points
which are both forward and backward recurrent.  By the Poincar\'e
Recurrence Theorem, $\mu((\Ld\cap V)\setminus A)=0$.

\begin{Definition}
A point $x\in A$ is {\em $\delta$-inaccessible} or simply {\em
inaccessible} if it is accumulated on both sides of $W^s_N(x)$ by
points in $\pi^s_x(A)$ and is accumulated on both sides of $W^u_N(x)$
by points in $\pi_x^u(A)$. Otherwise, $x$ is {\em accessible}.
\end{Definition}

Notice that this definition {\em does} depend on the choices
made. 

By a {\em rectangle} it is meant a Jordan curve made up of alternating
segments of stable and unstable manifolds, two of each. The segments
forming the boundary are its {\em sides} and the intersection points
of the sides are the {\em corners}. A rectangle is said to {\em
enclose} $p$ if it is the boundary of an open topological disk
containing~$p$.

\begin{Lemma}
\label{lem:accessible0}
Let $x\in A$ be an inaccessible point. Then there exist rectangles
enclosing $x$, having sides along the invariant manifolds of
hyperbolic periodic saddles in $V$ and having arbitrarily small
diameter.
\end{Lemma}
\begin{proof}
By assumption, there are points in $A$ whose stable manifolds are
arbitrarily near that of $x$ and to the right of it (think of stable
manifolds as roughly vertical and unstable manifolds as roughly
horizontal). Since points in $A$ are recurrent, by Katok's Closing
Lemma, it is possible to can find a periodic point, also to the right
of the stable manifold of $x$, whose stable manifold is
$N$-admissible. Proceeding like this, it is possible to find periodic
points, whose invariant manifolds are $N$-admissible, on all four
`quadrants' determined by the stable and unstable manifolds of $x$.
Clearly, segments of these manifolds form a rectangle enclosing $x$
and the choices can be made so that the rectangle has arbitrarily
small diameter.
\end{proof}

Notice it now follows from standard arguments that these periodic
saddles have transverse homoclinic intersections. In particular, the
closures of their invariant manifolds carry positive entropy. This
observation will be used below.

\begin{Lemma}
\label{lem:accessible1}
The set of accessible points in A has measure 0.
\end{Lemma}
\begin{proof}
Pick $x\in A$ and define $L^s\subset W^s_N(x)$ to be the set of endpoints
of the non-degenerate components of $W^s_N(x)\setminus\pi^s_x(A)$ and
$L^u\subset W^u_N(x)$ to be the set of endpoints of the non-degenerate
components of $W^u_N(x)\setminus\pi^u_x(A)$. Then both $L^s,L^u$ are
countable and it follows from the definitions that the set of
accessible points in $A$ lie on
$(\pi^s_x)\I(L^s)\cup(\pi^u_x)\I(L^u)$. Since each fiber of $\pi^s_x,
\pi^u_x$ contains at most one point which is both forward and backward
recurrent and since $\mu$ is non-atomic, it follows that the sets
$(\pi^s_x)\I(L^s),(\pi^u_x)\I(L^u)$, and thus the set of accessible
points, have measure~0.
\end{proof}

\begin{figure}
\lab{C}{C}{b}
\lab{p}{x}{l}
\lab{e}{\frac{\veps}{2}}{r}
\lab{s}{W^s_N(x)}{l}
\lab{u}{W^u_N(x)}{l}
\lab{z}{z_1}{}
\lab{w}{z_2}{}
\lab{y}{z_3}{}
\lab{C1}{C_1}{l}
\lab{C2}{C_2}{l}
\lab{b}{\beta}{}
\lab{enlarge}{\text{enlarge}}{b}
\pichere{0.5}{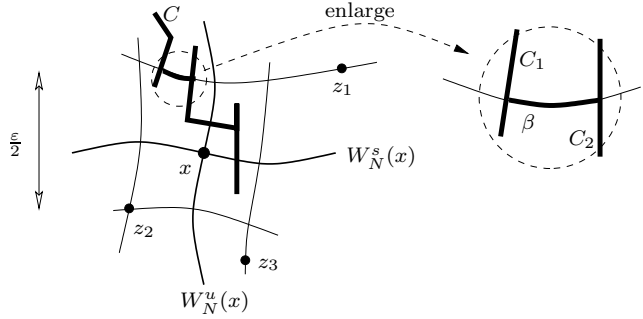}
\caption{An inaccessible point and a nearby continuum.}
\label{fig:inaccessible}
\end{figure}

In the proof below the following concept will be used.  Let $\gamma$ be
an open arc. A connected set $C$ is said to {\em cross} $\gamma$ if
there exist two disjoint connected subsets $C_1,C_2\subset
C\setminus\gamma$, on different 
sides\footnote{In order to make precise sense of the 
`different sides of $\gamma$' one has to invoke the Jordan Curve
Theorem.} of~$\gamma$, and a (possibly
degenerate) sub-arc $\beta\subset\gamma\cap C$, such that
$C_1\cup\beta\cup C_2$ is connected (see the diagram on the right in
Figure~\ref{fig:inaccessible}). 
\begin{Lemma} 
\label{lem:accessible2}
Suppose $x\in A$ is an inaccessible point. Given
$\veps>0$, there exists $0<r<\veps$ such that, if $z,w$ are any two points satisfying 
$d(x,z)<r$ and $d(x,w)>\veps$, then $z,w$ are  not zero-entropy equivalent. 
\end{Lemma}
\begin{proof}
By Lemma~\ref{lem:accessible0}, there exists a rectangle $R$ enclosing
$x$ with $\diam(R)<\veps$ and with sides along the invariant manifolds
of periodic points. Moreover, by the observation following the proof
of that lemma, these periodic points have transverse homoclinic
intersections and, therefore, their invariant manifolds carry
entropy. Pick $0<r<d(x,\partial R)$ (so that $r<\veps$ since $R$ is a Jordan curve enclosing $x$). If $z,w$ are points satisfying the hypotheses of the lemma and $C$ is a continuum containing $z,w$ then $C$ must cross the boundary of
$R$ (see Figure~\ref{fig:inaccessible} where $z_i,\, i=1,2,3$
are periodic points). The proof of the lemma now follows from
the simple observations below:
\begin{enumerate}[a)]
\item If a continuum $C$ crosses the stable manifold of a periodic saddle
$y$ of period $k$, then there is a compact subset $K\subset C$ such that
$d_H(f^{nk}(K),\ol{W^u(y)})\to 0$ as $n\to\infty$ (and an analogous
statement holds exchanging the roles of stable and unstable manifolds).
\item If $K,K'$ are compact sets and $d_H(f^n(K),f^n(K'))\to 0$ as
$n\to\infty$, then $h(f,K)=h(f,K')$.
\end{enumerate}
\end{proof}

\begin{Lemma}\label{lem:nequivalent}
If $C$ is a continuum with $h_\pm(C)>0$ then every continuum that is close   enough to $C$ in the Hausdorff metric contains inequivalent points. In 
particular, carrying positive entropy is an open condition in the space of continua.
\end{Lemma}
\begin{proof}

Assume that $h(f,C)>0$: an analogous argument can applied to $f\I$ in case $h(f\I,C)>0$. Let $\mu$ be the measure given by
Lemma~\ref{lem:basic}. Just as before, it follows that $\mu$ is
non-atomic and hyperbolic. Choose $0<\d<1$, $\Ld, p$ and $V$ as in the
discussion above so that $\mu(A)>0$, where
$A$ is the set of recurrent points in $\Ld\cap V$. By
Lemmas~\ref{lem:basic} and~\ref{lem:accessible1}, there exists an
inaccessible point $x\in A$ lying on a non-degenerate continuum
$E\in\omega(C)$. Let $\veps=\diam(E)/10>0$ and let $0<r<\veps$ be given by
Lemma~\ref{lem:accessible2}. Since $E\in\omega(C)$, there exist $\rho>0$ and $k\in\N$ such that, if $D$ is any continuum with $d_H(D,C)<\rho$, then 
$d_H(f^k(D),E)<r$. Since $x\in E$ this implies that there exist points $z,w\in f^k(D)$ such that $d(z,x)<r$ and $d(w,x)>\veps$. For the latter inequality, notice that there is $y\in E$ such that $d(x,y)\ge\diam(E)/2=5\veps$ and that there is $w\in f^k(D)$ such that $d(w,y)<r<\veps$. Thus, $d(w,x)\ge d(x,y)-d(w,y)>5\veps-\veps=4\veps$.
It now follows
from Lemma~\ref{lem:accessible2} that $0<h_\pm(f^k(D))=h_\pm(D)$.
\end{proof}

\medskip
\noindent
{\em Proof of Theorem~\ref{thm:musc}.}  Monotonicity is obvious from the definition of zero-entropy equivalence. Upper semi-continuity now follows easily from Lemma~\ref{lem:nequivalent}: if $x_n\to x$, $y_n\to y$ and $C_n$ are continua containing both $x_n,y_n$ with $h_\pm(C_n)=0$, by compactness there is a convergent subsequence $C_{n_k}\to C$, where $C$ is a continuum, $C\ni x,y$ and, by Lemma~\ref{lem:nequivalent}, $h_\pm(C)=0$. That the equivalence classes carry zero entropy also follows immediately from Lemma~\ref{lem:nequivalent}.\qed

\medskip\noindent
{\em Proof of Theorem~\ref{thm:main}.} Let $M_\sim$ denote the quotient of $M$ under the zero-entropy
equivalence. By the Moore-Roberts-Steenrod theorems, $M_\sim$ is a
generalized cactoid. Let $\pi\colon M\to M_\sim$ denote the quotient
map. Since the zero-entropy equivalence is clearly $f$-invariant, $f$
projects to a homeomorphism $f_\sim\colon M_\sim\barrow$. 
To finish the proof of the theorem all there is left is to show is that
$f_\sim$ is tight. The proof of Theorem 17 in~\cite{Bo} shows that for
any compact set $C\subset M_\sim$
\[h(f, \pi^{-1}(C))\le h(f_\sim,C)+\sup_{x\in C} h(f, \pi^{-1}(x)).\] 
Since the fibers of $\pi$ do not carry entropy, it follows that 
$h(f, \pi^{-1}(C))\le h(f_\sim,C)$. 
If $C\subset M_\sim$ is a non-degenerate continuum, then $\pi\I(C)$ is also a
continuum (since $\pi$ is monotone) and contains inequivalent
points (since $C$ is non-degenerate). Therefore $h(f, \pi^{-1}(C))>0$,
which completes the proof.
\qed

\section{Proof of Lemma~\ref{lem:basic}}
\label{sec:app}

The lemmas below, whose proofs can be found in~\cite{KH} (Lemma 4.5.2)
and~\cite{Ma} (Lemma 11.8) respectively, will be central to
what follows.

\begin{Lemma}
\label{lem:KH}
Let $X$ be a compact metric space, $f\colon X\barrow$ a homeomorphism,
$E_n$ $(n,\varepsilon)$-separated sets, 
\[\nu_n=\frac{1}{\card(E_n)}\sum_{p\in E_n}\delta_p\]
the uniform $\delta$-measures on $E_n$ and  
\[\mu_n=\frac{1}{n} \sum_{i=0}^{n-1}f_{\ast}^i\nu_n.\]
Then there exists an accumulation point $\mu$ of $\{\mu_n\}_{n\in\N}$
in the \weaks topology such that $\mu$ is an $f$-invariant Bpm
satisfying
\[h_{\mu}(f)\geq \limsup_{n \rightarrow \infty}
\frac{1}{n}\ln\card(E_n).\]\qed
\end{Lemma}

\begin{Lemma}[Pliss]
\label{Pliss}
For every $\lambda\in\R$, $\eta>0$, $H>0$ there are
$N_0=N_0(\lambda, \eta, H)\in\N$, $\delta=\delta(\lambda,\eta,H)>0$
such that if $a_1,\ldots,a_N$ are real numbers, $N\geq N_0$,
$|a_n|\leq H$, and 
\[\sum_{i=1}^{N}a_i \leq N\lambda\] 
then there are 
$1\leq n_1 < n_2 <\ldots < n_\ell \leq N$ such that 
\[\sum_{i=n_j +1}^{n}a_i \leq (n-n_j)(\lambda+ \eta)\]
for $j=1,\ldots,l$  $n_j<n \leq N$. Moreover $\ell/N \geq \delta$.\qed
\end{Lemma}

\smallskip
Recall that $\cK$ denotes the space of compact subsets of the compact
metric space $X$ with the Hausdorff metric.
A continuous map ({\em resp.} homeomorphism) $f\colon X\barrow$
induces a continuous map ({\em resp.} homeomorphism) $\hat
f\colon\cK\barrow$ in the obvious way: $\hf(K)=f(K)$. (It may seem
pedantic to differentiate between $f$ and $\hf$, but it is useful to
do so in order to avoid confusion below.) Recall that $D(n,\veps,K)$
denotes the minimum number of sets with $d_n$-diameter less than $\veps$
needed to cover $K$. Notice that, as a function of $K$ (with values in
$\N$), $D(n,\veps,K)$ is upper semi-continuous. The following useful
inequality is an easy consequence of the definition:
\begin{equation}
D(m+n,\veps,K)\le D(m,\veps,K)\cdot D(n,\veps,f^m(K)).
\label{subadd}
\end{equation}

Define three Borel measurable functions $h_{n,\veps}, h_\veps,h\colon\cK\to
[0,\infty]$ by 
\begin{eqnarray*}
h_{n,\veps}(K)&=&\frac{1}{n}\ln D(n,\veps,K) \\
h_\veps(K)&=&\lim_{n\to\infty}h_{n,\veps}(K)\\ 
h(K)&=&\lim_{\veps\to 0}h_\veps(K)
\end{eqnarray*}
Notice that $h(K)=h(f,K)$. That the limit as
$n\to\infty$ exists follows from inequality~(\ref{subadd}) and
that the limit as $\veps\to0$ exists follows from the monotonicity of
$h_{n,\veps}$ as a function of $\veps$.  Clearly, for any $K\in\cK$,
$h_{n,\veps}(K)\le h_{n,\veps}(X)$. As has already been observed, $h$ is
$\hf$-invariant, that is, for any compact set $K\in \cK$,
$h(K)=h(\hf(K))$.

Let $C \in \cK$ and define ${\mathfrak M}_{C}$ to be the set
of $\hf$-invariant Bpm's $\hat\mu$ on $\cK$ which are
\weaks limit points of the set of orbital averages of the form
\begin{equation}
\hat\mu_n= \frac{1}{n} \sum_{i=0}^{n-1} \delta_{\hf^i(C)}
\label{eq:orbitalmeasure}
\end{equation}
where $\delta_x$ denotes the Dirac mass at the point $x$. 
Observe that measures in $\fM_C$ have support contained in
$\omega(C)=\omega_{\hf}(C)$, the $\omega$-limit set of $C$ under~$\hf$.

\begin{Lemma}
\label{lem:basic0} 
Let $X$ be a compact metric space and $f\colon X\barrow$ be a homeomorphism and assume $C\in\cK$ satisfies $h(C)>0$. Then
$$\int_{\cK}h \; d\hmu>0$$
for every $\hmu\in{\fM}_{C}$
\end{Lemma}

\begin{proof}
Take $0<\eta<{h(C)}/{2}$ and set $\lambda= -h(C)+\eta$.
Choose
$\varepsilon_0>0$ and $N_1$ such that
\[ \ln D(n, \varepsilon_0,C)\geq -\lambda n>0\]
for every $n\geq N_1$.
Let $\Gamma_m$ denote the set of $K\in\cK$ such that $h_{m,\veps_0}(K) \geq
-\lambda -\eta$. Observe that, by upper semi-continuity, the
sets $\Gamma_m$ are closed, and therefore compact.

To make the notation less cumbersome, set $c_n=D(n,\veps_0,C)$. Set also
$c_0=1$ and
\[a_n=\ln\frac{c_{n-1}}{c_{n}}.\]
It follows from inequality~(\ref{subadd}) that there is $H=H(\veps_0)$,
independent of $n$, such that $|a_n| \leq H$.

Observe that
\[\sum_{i=1}^{N}a_i = \ln\frac{c_0}{c_N} = -\ln c_N \leq \lambda N.\]
for every $N \geq N_1$.
Choose $\delta=\delta(\lambda,\eta,\veps_0)>0$ and $N_0\geq N_1$ according
to Lemma~\ref{Pliss} for $\lambda$, $\eta$ and $H$ as above.
We are going to show that, for orbital average measures $\hmu_N$ defined by equality~(\ref{eq:orbitalmeasure}),
\(  \hmu_N(\Gamma_m)  \geq \frac{\delta}{2}\)
for all integers $m$ and $N\geq N_0$ such that ${m}/{N} <{\delta}/{2}$.

Take $1\leq n_1 < n_2 <\dots< n_l \leq N$
according to Lemma~\ref{Pliss} such that
\[\ln\frac{c_{n_j}}{c_n}=\sum_{i=n_j +1}^{n}a_i \leq (\lambda+\eta)(n-n_j)\]
for $j=1,\ldots, \ell$, $n_j<n \leq N$. Inequality~(\ref{subadd}) states
that $c_n\leq c_{n_j}\cdot D(n-n_j,\veps_0, f^{n_j}(C))$, from which it
follows \(h_{n, \veps_0}(\hf^{n_j}(C)) \geq -\lambda -\eta\), for $0
\leq n \leq N-n_j$.

Write $\ell=k_1+k_2$, where $k_1$ is the number of values $n_j$ such that
$0 \leq n_j \leq N -m$.
Then,
\[ \hmu_N(\Gamma_m) \geq  \frac{k_1}{N}.\]
On the other hand, by Lemma~\ref{Pliss},
\[\frac{k_1+ k_2}{N} \geq \delta \]
and it follows that
\[  \hmu_N(\Gamma_m)  \geq \delta - \frac{m}{N} \geq \frac{\delta}{2}.\]
It now follows from \weaks convergence and from the fact that
$\Gamma_m\subset\cK$ is compact that $\hmu(\Gamma_m)\ge\delta/2$, so
that
\[ \int h_{m,\veps} \; d\hmu \geq\int h_{m,\veps_0} \; d\hmu \geq (-\lambda 
-\eta)\frac{\delta}{2}.\]
As was pointed out above, $h_{m,\veps}$ is uniformly bounded by
$h_{m,\veps}(X)\to h_{\veps}(X) <\infty$. Applying the Dominated
Convergence Theorem as $m\to\infty$ and then the Monotone Convergence
Theorem as $\veps\to 0$, it follows that
\[ \int h \; d\hmu \geq (-\lambda -\eta)\frac{\delta}{2}>0\]
as was claimed.
\end{proof}

Denote by ${\mathcal F}$ the set of pairs $(p,C)$ such that $C\in\cK$
and $p \in C$. Clearly $\cF$ is a compact subset of $X\times\cK$.  A
homeomorphism $f\colon X\barrow$ induces a map $\tf\colon{\mathcal
F}\barrow$ setting $\tf(p,C)=(f(p),\hf(C))$. Let $\pi\colon{\mathcal
F}\to X$ and $\Pi\colon{\mathcal F}\to\cK$ denote the canonical
projections.  As usual with projections, the induced push-forward maps
$\pi_{\ast}$ and $\Pi_{\ast}$ send ergodic measures to ergodic
measures.

It is now possible to give a 

\medskip\noindent {\em Proof of Lemma~\ref{lem:basic}.}  Let
$\hlambda$ be an ergodic Bpm on $\cK$ such that $h >0$ almost everywhere with
respect to $\hlambda$ and $\supp\hlambda\subset\omega(C)$.  That such
a measure exists follows from the previous lemma and the Ergodic
Decomposition Theorem (Theorems II.6.1 and II.6.4 in~\cite{Ma}). From
these, it also follows that it is possible to choose a generic point
$A$ of the measure $\hlambda$ (that is, a compact $A\in\cK$ such that
$\hlambda$ is the \weaks limit of \((1/n) \sum_{i=0}^{n-1}
\delta_{\hf^i(A)}\)), with $h(A)>0$.

Choose $\veps>0$ and maximal $(n,\varepsilon)$-separated sets $E_n \subset A$ 
such that 
\[\limsup_{n\to\infty}\frac{1}{n}\ln\card(E_n)>0\]
and consider the measures
on $\cF$ given by
\[\tnu_n=\frac{1}{\card(E_n)} \sum_{x \in E_n}
\delta_{(x,A)}\quad\text{ and }\quad \tlambda_n= \frac{1}{n}
\sum_{i=0}^{n-1} \tf^i_\ast\nu_n.\] 
Setting $\lambda_n=\pi_\ast(\tlambda_n)$, Lemma~\ref{lem:KH} implies that there exists a \weaks convergent subsequence $\lambda_{n_k}\to\lambda$
so that $h_\lambda(f)>0$. Now, let $\tlambda$ be a \weaks limit of $\{\tlambda_{n_k}\}$. Notice that $\pi_\ast(\tlambda)=\lambda$ and $\Pi_\ast(\tlambda)=\hlambda$.
Let $N=\{x\in X\,\,|\,\, x\text{ belongs to a
non-degenerate set in } \omega(C)\}$. We claim that $\lambda(N)=1$. To
see why this is so, let $\hat N=\{K\in\cK\,\,|\,\, K\in\omega(C) \text{
is a non-degenerate set}\}$ and $\tilde N=\Pi\I(\hat
N)=\{(x,K)\in\cF\,\,|\,\, K\in\omega(C) \text{ is a non-degenerate
set}\}$.  Notice that $N=\pi(\tilde N)$. Since $h >0$
$\hlambda$-almost everywhere, $\hlambda(\hat N)=1$, so that
$\tlambda(\tilde N)=\tlambda(\Pi\I(\hat N)) =\Pi_\ast\tlambda(\hat
N)=\hlambda(\hat N)=1$.  On the other hand,
$\lambda(N)=\pi_\ast\tlambda(\pi(\tilde N)) =\tlambda(\pi\I(\pi(\tilde
N))) \ge\tlambda(\tilde N)=1$.

Applying the Ergodic Decomposition Theorem to $\lambda$ and using the
affine dependence of the measure theoretic entropy on the measure (see
formula S.2.5, p.\ 670 of~\cite{KH}) an ergodic Bpm $\mu$ on
$X$ with the desired properties is obtained.  \qed

\section{Proof of Theorem~\ref{thm:expansive}}
\label{sec:expansive}

Recall that $f$ is continuum-expansive if there is $\varepsilon >0$
such that for every non-degenerate continuum $C$ there is $n \in {\Z}$
such that $\diam f^n(C) \geq \varepsilon$.  In what follows this $\varepsilon$
is kept fixed.  The following convexity lemma will be
needed:
\begin{Lemma}
\label{convex}
There exists $r>0$ such that, for every continuum $C$ and every integer
$N>0$ satisfying $\diam(C)<r$ and $\diam(f^N(C))<r$, it follows that 
$\diam(f^n(C))<\varepsilon$ for every $0 \leq n \leq N$.
\end{Lemma}
\begin{proof}
Let $B_N(x,\varepsilon)$ denote the set of points $y$ such that $d(f^j(x),
f^j(y))\leq \varepsilon$ for every $0 \leq j \leq N$. If the lemma is
false, there are a sequence of continua $\Gamma_n$ and a sequence of
integers $N_n$ such that $\diam(\Gamma_n)<1/n$,
$\diam(f^{N_n}(\Gamma_n))<1/n$ and such that there exist $0<j_n<N_n$
for which $\diam (f^{j_n}(\Gamma_n))\geq \varepsilon$.
 
Choose points $x_n\in\Gamma_n$ and let $C_n$ be the connected
component of the set $B_{N_n}(x_n,\varepsilon /2)\cap \Gamma_n$
containing $x_n$.  Then, $\diam (f^{i}(C_n))\leq \varepsilon$ for
every $0\leq i\leq N_n$ and there is $0 \leq m_n \leq N_n$ such that
$\diam(f^{m_n}(C_n))\geq \varepsilon/2$.  Because $f$ is a
homeomorphism and $M$ is compact, it follows that $m_n\to\infty$ and
$N_n-m_n\to\infty$.

Set $D_n=f^{m_n}(C_n)$ and let $n_k \rightarrow \infty$ be a sequence
of integers such that $D_{n_k}$ converges in the Hausdorff topology to
the continuum $D$. Observe that $\diam(D)\geq \varepsilon/2$. On the
other hand, since $\diam(f^i(D_{n_k}))\leq \varepsilon$ for
$-m_{n_k}\leq i \leq N_{n_k}- m_{n_k}$, we get $\diam(f^i(D))\leq
\varepsilon$ for every $i \in {\Z}$, contradicting
continuum-expansivity.
\end{proof}

\medskip\noindent {\em Proof of Theorem ~\ref{thm:expansive}.}  Take
$5\rho<r<\varepsilon$ where $r$ and $\varepsilon$ are as in
Lemma~\ref{convex}. By continuum-expansivity and compactness there is
$N$ such that for every continuum $C$ such that $\diam(C)\ge\rho$
there is $|m|\leq N$ such that $\diam(f^m(C))\geq \varepsilon$.

Let $C$ be a continuum with $\diam(C)>\rho$ and assume there is $0\le
m\le N$ such that $\diam(f^m(C))\ge\veps$.  Lemma~\ref{convex} implies
$\diam(f^N(C))>r>5\rho$ so that it is possible to choose continua
$C_{1,1},C_{1,2}\subset f^N(C)$ such that $\rho <\diam(C_{1,i})<2\rho$
for $i=1,2$ and $d(C_{1,1},C_{1,2})>\rho$ (where
$d(A,B)=\inf\{d(a,b)\,\,|\,\, a\in A, b\in B\}$).  By induction on $k$ we
construct continua $C_{k,i}$, $1\leq i\leq 2^k$ such that $\rho
<\diam(C_{k,i})<2\rho$, $d(C_{k,2i-1},C_{k,2i})>\rho$,
$f^{-N}(C_{k,2i-1})\subset C_{k-1,i}$ and $f^{-N}(C_{k,2i})\subset
C_{k-1,i}$.  Assume that $C_{k,i}$ was constructed. We are going to
construct $C_{k+1,2i-1}$ and $C_{k+1,2i}$.  By continuum-expansivity
there is $|m|\leq N$ such that $\diam(f^m(C_{k,i}))>\varepsilon$.  It
follows that $m>0$ because otherwise $\diam(f^{-N}(C_{k,i}))>r>5\rho$ and 
$\diam(C_{k,i})>r>5\rho$ 
by Lemma \ref{convex} and this, in turn, contradicts the properties of
$C_{k,i}$.  Applying Lemma \ref{convex} again it follows that \(\diam(f^N(C_{k,i}))>r>5\rho\), which is enough to construct
$C_{k+1,2i-1}$ and $C_{k+1,2i}$. 

Let $T_k$ be a set obtained choosing exactly one point of each
$C_{k,i}$ and set $S_k= f^{-kN}(T_k)$.  Observe that the cardinality
of $S_k$ is $2^k$. By construction, if $x,y \in S_k$ are different,
there are $i\leq k$ and $j\leq 2^i$ such that $f^{iN}(x) \in C_{i,j}$
and $f^{iN}(y) \in C_{i,j+1}$ or vice versa. It follows that
$d(f^{iN}(x),f^{iN}(y))>\rho$, that is, $S_k$ is $(kN,\rho)$-separated.
Therefore, $h(f,C)>0$.

An analogous argument holds with $f\I$ in place of $f$ in case $-N\le
m\le 0$. It follows that $h_\pm(C)>0$ provided $\diam(C)\ge\rho$.
This completes the proof because if $C$ is an arbitrary non-degenerate continuum,
by continuum-expansivity there is $q$ such that
$\diam(f^q(C))\ge\veps>\rho$.
\qed 

\medskip
The converse of Theorem~\ref{thm:expansive} is not true in general. If
$f$ is the cactoidal map obtained by gluing countably many copies of
the tight horseshoe described in Section~\ref{sec:example} at their
fixed points, with diameters decreasing to 0, then $f$ is tight but not
continuum-expansive. We do not know whether the converse holds if
attention is restricted to homeomorphisms of manifolds, for example. 

Given $\varepsilon>0$, define the $\varepsilon$-length
$L_{\varepsilon}(\alpha)$ of a path $\alpha\colon [0,1] \rightarrow X$
as the smallest number $n$ for which it is possible to write $\alpha$
as the concatenation of $n$ paths
$\alpha=\alpha_1\cdot\ldots\cdot\alpha_n$, with
$\diam(\alpha_i)<\veps$. It is not hard to show that tight maps expand
the lengths of curves in the following sense: if $f\colon X\barrow $
is a tight homeomorphism then for every path $\alpha\colon [0,1]\to X$
there is $\varepsilon>0$ such that
\(\limsup_{n\rightarrow\infty}\frac{1}{n} \ln L_{\varepsilon}(f^n
\circ \alpha)>0\). If $f$ is continuum-expansive then the 
proof above shows that
$\varepsilon$ can be chosen uniformly for every $\alpha$.

\section{Problems for further research}
\label{sec:further}
There are several questions that arise naturally:
\begin{enumerate}[a)]
\item The Moore-Roberts-Steenrod theorems also give conditions under
which the quotient space of a surface under a monotone upper
semi-continuous decomposition is homeomorphic to the original
surface. (In the case of the sphere the condition is that the
complement of every decomposition element be connected.) It would be
interesting to give conditions on the map to ensure that this is the
case for the associated zero-entropy decomposition. 
\item In the horseshoe example, the measure of maximal entropy gives
rise to a Riemann surface structure together with an integrable
meromorphic quadratic differential on the quotient sphere (with one
essential singularity on the accumulation point of 1-pronged
singularities (the poles of the quadratic differential)).  The
quotient map is a Teichmuller mapping with respect to this
structure. We believe it is possible to show this is the case whenever
$f$ satisfies Axiom A and any two of its basic sets which do not
reduce to periodic orbits are unrelated by Smale's partial order
(that is, if $\Lambda,\Lambda'$ are any two distinct basic sets (which
are not periodic orbits), then $W^s(\Lambda)\cap
W^u(\Lambda')=\emptyset$).  This is related to a result of
Bonatti-Jeandenans~\cite{BL}.
\item In work in preparation, the first author shows that, under
certain assumptions, it is possible to take limits of sequences of
generalized pseudo-Anosovs. The maps obtained are tight and the
Riemann surface structure is still present, but the quadratic
differentials may be only $L^1$ and have dense sets of `zeros' and
`poles.' To find conditions under which the zero-entropy quotient of a
surface diffeomorphism is the limit of generalized pseudo-Anosov seems
an interesting problem. If this is the case, the quotient map should
be thought of as a piecewise affine model for the diffeomorphism. 
\end{enumerate}

\end{document}